\documentclass[11pt]{article}
\usepackage{amssymb}    
\usepackage{amsmath}    
\usepackage{amsthm} 
\usepackage{mathtools}
\usepackage{color}  
\newtheorem{theorem}{Theorem}

\newtheorem{lemma}{Lemma}
\newtheorem{definition}{Definition}

\def\D{{\cal D}}

\def\R{{\mathbb R}}
\def\Rd{{\mathbb R}^d}

\def\HH{{\cal H}}
\def\PROB {{\mathbb P}}
\def\EXP {{\mathbb E}}
\def\IND{{\mathbb I}}
\def\Var{{\mathbb Var}}

\def\sgn{\mathop{\rm sgn}}

\allowdisplaybreaks

\begin{document}
\begin{titlepage}
\thispagestyle{empty}
\setcounter{page}{0}

\title{Distribution-free tests for lossless feature selection in classification and regression}
\author{
L\'aszl\'o Gy\"orfi\thanks{Department of Computer Science and Information Theory, Budapest University of Technology and Economics, gyorfi@cs.bme.hu. The research of L\'aszl\'o Gy\"orfi was supported by the National Research, Development and Innovation Fund of Hungary 
under the 2023-1.1.1-PIACI-F\'OKUSZ-2024-00051 funding scheme.}
\and
Tam\'as Linder\thanks{Department of Mathematics and Statistics, Queen's University, Kingston, ON Canada, tamas.linder@queensu.ca. Tam\'as Linder's research was supported in part by the Natural Sciences and Engineering Research Council (NSERC) of Canada.}
\and
Harro Walk\thanks{Institut f\"ur Stochastik und Anwendungen, Universit\"at Stuttgart, harro.walk@mathematik.uni-stuttgart.de.}}
\maketitle

\begin{abstract} 
We study the problem of lossless feature selection for
a $d$-dimensional feature vector $X=(X^{(1)},\dots ,X^{(d)})$ and label $Y$ for binary classification as well as nonparametric regression. 
For an index set  $S\subset \{1,\dots ,d\}$, consider the selected $|S|$-dimensional feature subvector $X_S=(X^{(i)}, i\in S)$. If $L^*$ and $L^*(S)$ stand for the minimum risk based on $X$ and $X_S$, respectively, then $X_S$ is called lossless if $L^*=L^*(S)$. For classification, the minimum risk is the Bayes error probability, while in  regression,  the minimum risk is the residual variance. 
We introduce nearest-neighbor based test statistics   to test the hypothesis that $X_S$ is lossless. This test statistic is an estimate of the excess risk $L^*(S)-L^*$.
Surprisingly,  estimating this excess risk turns out to be
a functional estimation problem that does not suffer from the curse of dimensionality 
in the sense that the convergence rate does not depend on the dimension  $d$.
For the  threshold $a_n=\log n/\sqrt{n}$, the  corresponding tests are proved to be consistent under conditions on the distribution of $(X,Y)$ that are significantly milder than in previous work. Also, our threshold is universal (dimension independent),  in contrast  to earlier methods where for large $d$ the threshold becomes too large to be useful in practice.

\end{abstract}

\noindent


\noindent


\thispagestyle{empty} 

\end{titlepage}


\section{Introduction}

In this paper we study the problem of lossless feature selection for classification and nonparametric regression.

Binary classification deals with the problem of deciding on a $\pm 1$-valued random label $Y$  based on a 
random feature vector $X$ taking values in $\Rd$, 
so that the risk, measured by  the decision error probability, is as small as possible. If the joint distribution of $X$ and $Y$ is known, then the optimal decision with minimum risk (error probability), called the
Bayes decision, can be derived. In the standard setup of classification, the joint distribution is unknown, but instead an observed random sample $(X_1,Y_1),\ldots,(X_n,Y_n)$ of $n$ independent copies of $(X,Y)$ is available from which an estimate of the Bayes decision is to be constructed. Although estimates (classification rules)  exist whose error probabilities converge to the optimum as $n\to \infty$ without any condition on the distribution of $X$ and $Y$, 
the convergence rate of the error probability of any such  classification rule to the Bayes error is very sensitive to the dimension of the feature vector.
This suggests that dimension reduction, also called feature selection,  is crucial before constructing a  classification rule. 

For nonparametric regression, $Y$ is a real-valued random variable with $\EXP[Y^2]<\infty$, 
the risk is the mean squared error,  and the minimum risk is the residual variance.

In this paper we are interested in feature selection that is lossless, i.e., does not incur any loss of information. 
For lossless feature selection, the minimum risk based  on the feature vector obtained by leaving out some components of $X$ and that  based on the original feature vector $X$, are equal.

To make this notion precise, let $S$ denote a proper subset of $\{1,\dots ,d\}$, and for an $\R^d$-valued  feature vector
\[
X=(X^{(1)},\dots ,X^{(d)}),
\]
consider the $\R^{|S|}$-valued  subvector picked out by $S$, given by 

\begin{align}
\label{XA}
X_S
=(X^{(i)}, i\in S).
\end{align}

Let $L^*$ and $L^*(S)$ denote the minimum risk  based on $X$ and $X_S$, respectively. Then   $X_S$ is called lossless if $L^*=L^*(S)$. One  goal in this paper is to construct a nonparametric (distribution free) test for the null hypothesis $\mathcal{H}_0 : L^*=L^*(S)$, i.e., for the null hypothesis  that the minimum risks based on $X$ and $X_S$ are equal. 
In our setup, the alternative hypothesis  $\mathcal{H}_1$ is that $L^*<L^*(S)$, i.e., $\mathcal{H}_1$, is the complement of the null hypothesis, and 
therefore there is no separation gap between the hypotheses.
It is not at all obvious that consistency is possible without any separation gap. 
Dembo and Peres \cite{DePe94} and Nobel \cite{Nob06} characterized  hypotheses pairs that admit strongly  consistent tests, i.e.,  tests that, with probability one, only make finitely many Type I and II errors. This property is called discernibility.
As an illustration of the intricate nature of the discernibility concept,
Dembo and Peres \cite{DePe94} demonstrated an exotic example, where the null hypothesis is that the mean of a random variable is rational, while the alternative hypothesis is that this mean minus $\sqrt{2}$ is rational.
(See also  Cover \cite{Cov73} and Kulkarni and Zeitouni \cite{KuZe91}.)
The discernibility property shows up in
Biau and Gy\"orfi \cite{BiGy05} (testing homogeneity),
Devroye and Lugosi \cite{DeLu02} (classification of densities),
Gretton and Gy\"orfi \cite{GrGy10} (testing independence),
Gy\"orfi and Walk \cite{GyWa12} and Gy\"orfi, Linder and Walk \cite{GyLiWa23} (testing conditional independence),
Morvai and Weiss \cite{MoWe21} and Nobel \cite{Nob06} (classification of stationary processes),
among others.

Consistent tests can be constructed by estimating the functionals $L^*$ and $L^*(S)$ using distribution-free consistent nonparametric estimates, i.e., estimates that (in some sense) converge to the target functional as $n\to \infty$. For example,  Gy\"orfi and Walk  \cite{GyWa17}  consider  the Bayes error probability, Devroye et al.\ \cite{DeGyLuWa18} the residual variance, Berrett et al.\ \cite{BeSaYu17}, Delattre and Fournier \cite{DeFo18}, and Kozachenko and Leonenko \cite{KoLe87} the differential entropy, Gretton and Gy\"orfi \cite{GrGy10},  Silva and Narayanan \cite{SiNa10}, and Wang et al.\ \cite{WaKuVe05} the mutual information, and Beirlant et al.\ \cite{BeDeGyVa01} the total variation.

Most of these estimators are based on consistent nonparametric estimators of the corresponding function (e.g., for the residual variance, the estimator of the regression function). However, the rate of convergence of such functional estimators is determined by the rate of convergence of the corresponding function estimator,  which can be slow. In general, estimating the function itself is a harder problem than estimating the corresponding functional because the variance of the functional estimate can be of order $O(1/n)$, as is the case in this work, see for example inequalities  (\ref{8'}), (\ref{V1}) and (\ref{V2}) in our analysis. 
In fact, it surprisingly  turns out that  functional estimators with good rate of convergence  are based on non-consistent function estimates, e.g., the 1-NN based estimators of residual variance \cite{DeGyLuWa18} and differential entropy \cite{KoLe87}. 

In this paper, we introduce  a distribution-free test which uses an estimate of the difference of the minimum risks $L^*-L^*(S)$, called the excess risk. 
This is in contrast to estimating $L^*$ and $L^*(S)$ separately and then taking the difference of these two estimates.  This is an important step since all separate estimates of  $L^*$  and   $L^*(S)$ suffer from the curse of dimensionality.
We propose a nearest-neighbor based test statistic $T_n$  and the threshold  $a_n= \log n/\sqrt{n}$. 
The null hypothesis $L^*=L^*(S)$ is accepted if $T_n\le a_n$,  and  otherwise rejected.
Our main results shows that the resulting tests (for classification and regression, respectively) are consistent in the sense that  the Type I and II errors converge to zero as the sample size $n$ tends to infinity. 
Surprisingly, estimating the excess risk turns out to be a  functional estimation problem where there is no curse of dimensionality in the sense that the convergence rate  (under the null hypothesis) depends on  $d$ only  through  a multiplicative factor (see  \eqref{eq-nod} at the end of the proof of Theorem~\ref{expo}.) 

For the alternative hypothesis, the analysis is relatively straightforward because the consistency of the test naturally follows from the consistency of the minimum risk estimates. 
However, under the null hypothesis, the problem is more challenging.  On the one hand, the estimation error $T_n-\EXP[\, T_n\, ]$ is of order $1/\sqrt{n}$. On the other hand, the absolute value of the bias $\EXP[\, T_n\, ]$ can be of order $n^{-1/d}$ ($d$ is the dimension of the feature vector $X$), which can be  much larger than  $1/\sqrt{n}$. The resulting threshold is about $a_n= \log n /n^{1/d}$, which  becomes impractically large for high-dimensional $X$, because the test rarely rejects the null hypothesis, see, e.g.,  Gy\"orfi and Walk  \cite{GyWa17}  (classification) and Devroye et al.  \cite[Section 3]{DeGyLuWa18}  (regression). To address this problem, we  develop  statistics $T_n$  such that the bias of the estimator of $L^*(S)$ is much larger than the bias of the estimator of $L^*$, resulting in a negative bias $\EXP[\, T_n\, ]$. This approach allows us to set a universal, dimension independent  threshold $a_n$, yielding smaller and more practically useful threshold values.

We note that if for a threshold sequence  $a'_n$, $T_n < a'
_n$ is a
consistent test, then so is $f_n(T_n) < f_n(a'_n)$ for any sequence of increasing scaling 
functions $f_n$, so by replacing the test statistic $T_n$ by $f_n(T_n)$, the
threshold $a_n=f_n(a'_n)$ can easily be made dimension independent. 
However, in this case $f_n(T_n)$ would not be 
a consistent estimate of the excess risk $L^*-L^*(S)$, which would make a test based on $f_n(T_n)$  much less attractive.

Under mild smoothness conditions, the main results in this paper show the consistency of the constructed tests for the lossless feature selection property for classification and nonparametric regression.
However, one may in addition be interested in the power of these tests. 
Without any separation gap between the hypotheses, for any test the rate of convergence of the Type II error (power) can be arbitrarily slow. In future work one may consider a formulation  of  separation that makes it possible to derive sharper bounds on the power on these tests. For example, with the null hypothesis $\mathcal{H}_0$ as above, one may consider the alternative hypothesis  $\mathcal{H}_1$ to be $L^*+\delta <L^*(S)$ for some fixed but unknown $\delta>0$ and investigate the power of the test in this setup as a function of $n$ and $d$.

The paper is organized as follows. In Section~\ref{sec-class} we introduce a novel $k$-nearest-neighbor based test statistic and threshold for  lossless feature selection in  classification and state the consistency of the corresponding test under a mild Lipschitz-type condition (Theorem~\ref{expo}). We also state a lemma (Lemma~\ref{MAD})  that gives sharp upper and lower bounds on the first absolute moment of the average of $n$ independent and identically distributed $\pm1$-valued random variables, which plays  an important  role in proving the test's consistency under the null hypothesis. In 
Section~\ref{sec-regr} we introduce a similar test for  lossless feature selection for nonparametric regression and state its consistency (Theorem~\ref{expo*}). Finally, the proofs are presented in Section~\ref{sec-proofs}.

\section{Lossless feature selection for classification}

\label{sec-class}

The task of binary classification is to decide on the $\pm 1$-valued random variable $Y$ given an $\R^d$-valued random vector  $X$ by finding a decision function $g$, defined on the range of $X$, such that $g(X)=Y$ with large probability.
If $g: \R^d \to \{-1,1\}$  is an arbitrary measurable decision function, then its error probability is denoted by
\[
L(g)=\PROB\{g(X)\ne Y\}.
\]

Let 
\[
m(x)=\EXP[\, Y \mid X=x].
\]
It is well-known that the so-called Bayes decision $g^*$, the decision function that  minimizes the error probability $L(g)$,  is given by
\[
g^*(x) =\sgn m(x),
\]
where  $\sgn(z)= 1$ if $z\ge 
0$ and $\sgn(z)=-1$ if $z<0$ for any 
$z\in \R$.

The minimum error probability, also called the Bayes error probability, is given by 
\[
L^*=\PROB\{g^*(X)\ne Y\}=\min_g L(g), 
\]
which can be considered as the minimum risk in classification.

Constructing the Bayes decision requires the knowledge of the  distribution of $(X,Y)$. Typically, this distribution is unknown and instead one observes the training samples
\[
\D_n=\{(X_1,Y_1), \dots ,(X_n,Y_n)\}, 
\]
consisting of independent and identically distributed (i.i.d$.$) copies of $(X,Y)$, which arrive in a stream with increasing sample size $n$. The monograph by Devroye et al. \cite{DeGyLu96} provides constructions of classification algorithms, based on the data $\D_n$,  that are universally consistent in the sense  that the error probability of these algorithms tends to the Bayes error probability for all distributions of $(X,Y)$ as $n\to \infty$.
However, the rate of convergence of the error probabilities heavily depends on regularity (smoothness) properties of the function $m$ and on the dimension $d$ of the feature vector $X$.
Detecting a subset of ineffective features that, in the presence of the other features, has no influence on $L^*$  allows lossless reduction of the dimension.

This section deals with testing the  hypothesis of  ineffectiveness of  specific features.
The test uses an estimate of the difference of the Bayes error probabilities with and without these features.

As before, for  $S\subset \{1,\ldots,d\}$ let $L^*(S)$ denote the Bayes error probability when $Y$ is estimated using the subvector $X_S=(X_i, i\in S)$. Our aim is to construct a test for the hypothesis $L^*=L^*(S)$, i.e., the hypothesis that  the  subset  of neglected features  $(X^{(i)}, i\notin S) $ do not provide more information beyond what is contained in $X_S$.
A set $S$ with this property provides  lossless feature selection. Note, however, that the neglected features  $(X^{(i)}, i\notin S) $ may still be informative, e.g., it is possible that $L^*=L^*(S)= L^*(S^c)$.

Most dimension reduction algorithms for classification are modified versions of principal component analysis (PCA), where one looks for a linear transformation of the feature vector into a lower dimensional subspace. In the resulting optimization problem the Bayes error probability is replaced by a smooth error proxy, see, e.g.,  Siblini et al. \cite{SiKuMe21} and Tang et al.\cite{TaAlLi14}.

A particular way of dimension reduction is  feature selection. 
Similarly to the PCA, in feature selection algorithms, instead of the Bayes error probability, one usually considers other, more treatable criteria, see  Guyon and Elisseeff \cite{GuEl03}.
Usually, such  feature selection algorithms are looking for individual relevant features.
Another, more direct goal is to find a  good feature set $S$ that has small size $|S|$ and  error probability  $L^*(S)$ that is close to $L^*$.
If one fixes $|S|=d'$ for some integer $1\le d' <d$, the problem is to find a $d'$-element selection $S$ with minimal Bayes error $L^*(S)$. 
In this respect, the examples of Cover and Van Campenhout \cite{CoVa77} and Toussaint \cite{Tou71} show that 
the  $d'$ features that are individually the best do not necessarily constitute the best $d'$-dimensional vector, and 
therefore, every algorithm has to search exhaustively through all $d'$-element subsets $S$, see Chapter 32 in  \cite{DeGyLu96}, This procedure is prohibitively complex and instead one may search 
 for the smallest feature set $S$ that provides lossless feature selection, a task that may be easier. In this paper we are concentrating on the problem of testing whether or not a given $S$ provides lossless feature selection, but  we do not deal with the problem of searching for such a smallest feature set $S$. Finding effective algorithms for identifying such $S$ may be the subject of future work.

When searching for such an $S$ based on the training samples $\D_n$, the classification null-hypothesis $\HH_0$, defined as 
\begin{align}
\label{testp}
L^*(S)=L^*
\end{align}
must be tested. The null-hypothesis (\ref{testp})  means that the subset of neglected features  $(X^{(i)}, i\notin S) $ of the vector $X$ carries no additional information,
i.e., it has no additional predictive power.

An obvious approach to this problem is  to estimate $L^*(S)$ and $L^* $ from the data and to accept the hypothesis (\ref{testp}) if the difference of the estimates is small.
Unfortunately, it seems that no such estimates with fast rate of convergence are avaiable in the literature.
Antos et al.\ \cite{AnDeGy99}
 proved that without any regularity conditions, the rate of convergence  for any estimate of $L^* $ can be arbitrarily slow.
In view of this, our approach will be to  estimate directly the difference $L^*(S)-L^*$. This task will prove easier in the sense that we can construct  an estimate $T_n$ of $L^*(S)-L^*$, for which, 
without any condition, 
\begin{align}
\label{TTn}
T_n
&\to
L^*(S)-L^*  
\end{align}
in $L_1$, and if $L^*(S)-L^*=0$, then
\begin{align*}
\EXP[\,|T_n|\, ]
&\to
0
\end{align*}
with a nontrivial rate of convergence.
More specifically, we introduce an estimate $T_n$ and a threshold $a_n\to 0$ such that (\ref{TTn}) holds and if $L^*(S)-L^*=0$, then
\begin{align*}
\lim_{n\to \infty} \PROB\{T_n>a_n\}
&=0.
\end{align*}

In order to detect ineffective features,  \cite{GyWa17} proposed nearest neighbor and partitioning based statistics and, for the threshold $a_n\approx 1/\sqrt{n}+n^{-2/(2+d)}$,  proved the consistency of the resulting tests under fairly restrictive conditions on the distribution of the pair $(X,Y)$.   For the case of large $d$, which is the focus in  feature selection, the corresponding  $a_n$ is too large to be useful in practice.
In this section we weaken the restrictive conditions and significantly decrease the threshold value.

Let $S\subset \{1,\ldots,d\}$ be fixed  and introduce the notation 
\[
\widehat X
=X_S,
\]
\[
\widehat m(\widehat X)=\EXP[\,Y\mid \widehat X]=\EXP[\, m(X)\mid \widehat X]
\]
and
\[
\widehat L^*= L^*(S)=\PROB\{\sgn\widehat m(\widehat X)\ne Y\}.
\]
For any measurable $g:
\R^d \to \{-1,1\}$ we have
\begin{align}
L(g)-L^*
&=\EXP\left[\IND_{\{ g( X)\ne g^*(X)\}}|m(X)|\,\right],
\label{g}
\end{align}
where $\IND$ denotes the indicator function
(see \cite[Theorem 2.2]{DeGyLu96}). 
Letting $g=-g^*$ in (\ref{g}), we obtain  $(1-L^*)-L^*=\EXP[\,|m(X)|\, ]$.
Therefore
\[
L^*
=\frac 12 \big(1-\EXP\big[|m(X)|\big]\big),
\]
and similarly
\[
\widehat L^*
=\frac 12 \big(1-\EXP\big[|\widehat m(\widehat X)|\big]\big),
\]
implying that 
\begin{align}
\widehat L^*-L^*
&=\frac 12 \left(\EXP\big[|m(X)|\big] -\EXP\big[|\widehat m(\widehat X)|\big] \right),
\label{l2}
\end{align}
Therefore the classification null hypothesis (\ref{testp}) is equivalent to 
\begin{align}
\label{testl}
\EXP\big[ |m(X)|\big] -\EXP\big[|\widehat m(\widehat X)|\big] =0.
\end{align}

To test the null hypothesis, we propose a nearest-neighbor-based test statistic.
Fix $x  \in \Rd$ and reorder the data $(X_1,Y_1),\ldots,(X_n,Y_n)$
according to increasing values of $\|X_i-x\|$, where $\|\cdot \|$ denotes the Euclidean norm on $\R^d$.
The reordered data sequence is denoted by
\begin{align}
\label{NN}
(X_{(n,1)}(x),Y_{(n,1)}(x)),\ldots,(X_{(n,n)}(x),Y_{(n,n)}(x)),
\end{align}
so that $X_{(n,k)}(x)$ is the $k$th nearest-neighbor of $x$.
In case of a tie, i.e., if $X_i$ and $X_j$ are equidistant
from $x$, then $X_i$ is declared ``closer'' if $i < j$.
In this paper we assume that ties occur with probability $0$.
This assumption can be enforced by endowing $X$ with a $(d+1)$th component $Z$ that  is independent of $(X,Y)$ and is uniformly distributed on $[0,1]$, see Section 11.2 in \cite{DeGyLu96}. 
In this situation, the training samples  are similarly modified.  We note that this procedure does not change either $L^*$ or $\widehat{L}^*$ and $m$ and $\widehat{m}$ remain unchanged as well. 

Choose an integer $k_n$ less than $n$.
The $k$-nearest-neighbor ($k$-NN) regression estimate of $m$ is
\begin{align}
\label{est}
m_{n}(x) =\frac{1}{k_n}\sum_{j=1}^{k_n} Y_{(n,j)}(x).
\end{align}
Analogously to (\ref{est}),  from the training subsamples
\[
\widehat \D_n=\{(\widehat X_1,Y_1), \dots ,(\widehat X_n,Y_n)\}
\]
(where $\widehat{X}_i$ is the subvector picked out from $X_i$ by the selected features $S$), we introduce the $k$-NN estimate $\widehat m_{n} $ of $\widehat m$:
\begin{align}
\label{est-1}
\widehat m_{n}(\widehat x)
&= \frac{1}{k_n}\sum_{j=1}^{k_n} \widehat Y_{(n,j)}(\widehat x). 
\end{align}

The $k$-nearest-neighbor  classification rule  is defined by the plug-in estimator 
\begin{align*}
g_{n}(x)=\sgn (m_n(x)).
\end{align*}
Assuming ties occur with probability $0$, if $k_n$ is chosen so that $k_n\to \infty$ and $k_n/n\to 0$, then without any other condition on the distribution of $(X,Y)$, the $k$-NN classification rule is strongly consistent, i.e., 
\begin{align*}
L(g_{n}) 
\to L^* \quad\text{a.s.},
\end{align*}
where $L(g_n)= \PROB\{g_n(X)\neq Y \mid \D_n\}$, see Sections 11.1 and 11.2 in \cite{DeGyLu96}.

To construct the test we assume for convenience that $2n$ i.i.d.\ training samples are available, so that  in addition to $\D_n=\{(X_1,Y_1), \dots ,(X_n,Y_n)\}$, the i.i.d.\ copies 
\[
\D'_n=\{(X'_1,Y'_1), \dots ,(X'_n,Y'_n)\}
\]
of $(X,Y)$ are also available ($(X,Y)$, $\D_n$, and $\D'_n$ are independent). Recalling \eqref{l2}, we estimate  the difference 
\begin{align*}
\EXP\Big[ |m(X)|-|\widehat m(\widehat X)|\Big] 
\end{align*}
by means of the  nearest-neighbor based test statistic
\begin{align*}
T_n
&=
\frac 1n \sum_{i=1}^n \left( Y'_i\sgn(m_n(X'_i))-|\widehat m_n(\widehat X'_i)|\right)
\end{align*}
with $k_n$ chosen as 
\[ 
k_n
=\lfloor \sqrt{\log n}\rfloor.
\]

One accepts the classification null hypothesis if
\begin{align*}
T_n\le a_n,
\end{align*}
where
\begin{align}
\label{an}
a_n
&=
\frac{\log n}{\sqrt{n}}.
\end{align}

To prove the consistency of this test, we need the following modified Lipschitz condition. It is a combined smoothness and tail condition that weakens two rather restrictive conditions that are used in the literature:  the Lipschitz continuity of $m$ and the condition that $X$ is bounded.  

\begin{definition}[Chaudhuri and Dasgupta \cite{ChDa14}, D\"oring et al. \cite{DoGyWa18}]
If $\mu$ stands for the distribution of $X$, then   $m$ satisfies the modified Lipschitz condition if there exists $C^*>0$ such that for any
$x,z\in \Rd$,
\begin{align*}
|m(x)-m(z)|\le C^*\mu(S_{x,\|x-z\|})^{1/d},
\end{align*}
where  $S_{x,r}=\{y\in \R^d : \|y-x\|\le r\}$ denotes  the closed Euclidean ball centered at $x$ having  radius $r$.
\end{definition}

The following theorem, which is one of the main results in this paper,  states the consistency of our test.

\begin{theorem}
\label{expo}
Let $d\ge 2$ and assume that ties occur with probability $0$. Then 
\begin{itemize}
\item[\rm (a)]
Under the classification alternative hypothesis,  one has
\begin{align*}
\lim_{n\to \infty}\PROB\{ T_n\le a_n\}
&=0.
\end{align*}
\item[\rm (b)]
If $\widehat m$  satisfies the modified Lipschitz condition and
the residual variance satisfies $\EXP\big[ (Y-\widehat m(\widehat X))^2\big] =1-\EXP[ \widehat m(\widehat X)^2] >0$, then
under the classification null hypothesis, one has
\begin{align*}
\lim_{n\to \infty}\PROB\{ T_n> a_n\}
&=0.
\end{align*}
\end{itemize}
\end{theorem}

\noindent\emph{Remarks:}
\begin{itemize}
\item[(i)] The proof of the theorem also works if instead of the modified
\mbox{Lipschitz} condition one assumes the standard \mbox{Lipschitz} condition and the boundedness of $X$. 

\item[(ii)] In the proof we make use of the  Efron-Stein inequality  \cite[Theorem 3.1]{BuLuMa13}. At the price of a  larger threshold $(\log n)^2/\sqrt{n}$, the Efron-Stein inequality can be replaced by an exponential concentration inequality, such as McDiarmid's  inequality \cite[Theorem 6.2]{BuLuMa13},  to obtain  strong consistency. (A test is called strongly consistent if, almost surely, with increasing sample size  it makes finitely many errors.) 
\end{itemize}

In the analysis we show that, under the null hypothesis, the expectation of the test statistic is negative, and in the proof we need to lower bound the expected $L_1$ norm of the corresponding regression estimate.
In contrast to the $L_2$ setup, where the difference of the second moment and the squared expectation is equal to the variance, here the difference of the first absolute moment and the absolute value of the expectation is much less than the standard deviation if the expectation is non zero.
In this respect the lower bound in the following lemma plays a crucial  role. The proof of the lemma is given in Section~\ref{expo-proof}. 

\begin{lemma}
\label{MAD}
Let $Z_1,\dots , Z_n$ be $\pm 1$-valued i.i.d.\ random variables with mean $a$ and variance 
$\sigma^2=1-a^2$.
Then,
\begin{align}
\label{C4}
|a|
+
\frac{\sqrt{2}}{n^{3/2}}\sigma^{n}
&\le 
\EXP\Big|\frac{1}{n}\sum_{i=1}^n Z_i \Big|
\le 
|a|
+
\frac{c^*}{n^{3/2}}\sigma^{n}
\end{align}
with $c^*=c^*(a)<\infty$ if $a\ne 0$.
\end{lemma}

We note that for $a=0$ the upper bound in \eqref{C4} is not valid, but the lower bound is. In fact, for $a=0$  the Berry-Esseen theorem  implies the asymptotics  
\begin{align*}
\EXP\left|\frac{1}{n}\sum_{i=1}^n Z_i \right|
&=
\sqrt{\frac{2}{\pi n}}+ O\left( \frac 1n\right).
\end{align*}

Lemma \ref{MAD} can be considered as a binomial distribution analogue of the second part of Lemma 5.8 
in Devroye and Gy\"orfi  \cite{DeGy85}, which deals with the normal distribution.
However, an application of that result,  together with the normal approximation to the binomial distribution using the Berry-Esseen theorem  according to the first half of this Lemma 5.8, does not give the bound of Lemma~\ref{MAD}  because
the normal approximation is too rough for our setting.

\section{Lossless feature selection for nonparametric regression}

\label{sec-regr}

An analogous problem can be posed in the context of nonparametric regression,
where $Y$ is a real-valued random variable with $\EXP[Y^2]<\infty$.
Here  the functional
\[
L^*=\EXP  \left[( Y-m(X) )^{2} \right].
\]
with $m(x)= \EXP[ Y\,|\, X=x]$, plays the central role. 
It is known that for any  measurable function $g:\R^d\to\R$, 
\[
\EXP \left[ ( Y- g(X))^{2}\right] = L^*  + \EXP  \left[( m(X) - g(X) )^{2} \right]
\]
and therefore 
\[
L^*= \min_{g} \EXP  \left[ ( Y-g(X) )^{2} \right],
\]
where the minimum is taken over all measurable functions $g:\R^d\to\R$. The functional $L^*$ is often  referred to as the residual variance; it is the minimum  mean squared
error in predicting $Y$ based on the observation  $X$.

For $S\subset \{1,\ldots,d\}$, the predictive power  of a subvector
$X_S=(X^{(i)}, i\in S)$ of $X$ is measured by the residual variance
\[
L^*(S) = \EXP \left[ (Y - \EXP[Y\mid X_S])^2 \right]
\]
that can be achieved using the features as explanatory variables.
For possible dimensionality reduction, one needs, in general, to test the regression null hypothesis $\HH_0$  that the two residual variances are equal:
\begin{equation}
\label{rtest1}
L^*=L^*(S).
\end{equation}
Again, a set $S$ with this property provides lossless feature selection.

Identifying a set of features $S$ with the property \eqref{rtest1} is equivalent to finding a set of features $\bar{S}=S^c$ that are irrelevant in inference.  Lei and  Wasserman \cite{LeWa14} introduced  the Leave-One-Covariate-Out (LOCO) value  
\begin{align*}
\mbox{LOCO}(S^c)
&=
L^*(S)-L^*,
\end{align*}
for characterizing the importance of a feature or a subset of features,  see also  Gan et al.\  \cite{GaZhAl23}, Verdinelli and  Wasserman \cite{VeWa23}, and Williamson et al.~\cite{WiGiCaSi21}.
These papers mostly deal with finding a single element set $S^c$ that minimizes $\mbox{LOCO}(S^c)$. Testing the null hypothesis (\ref{rtest1}) is equivalent to  testing $\mbox{LOCO}(S^c)=0$ a task of which \cite{WiGiCaSi21} says 
``Developing valid inference under this particular null hypothesis appears very difficult.''

As in De Brabanteret al.\  \cite{BrFeGy14},
a natural way of approaching this testing problem is by estimating both residual variances
$L^*$ and $L^*(S)$, and accept the regression null hypothesis
if the two estimates are close to each other.

For fixed $S\subset \{1,\ldots,d\}$, let $\widehat{X}=X_S$ and $\widehat{m}(\widehat{x})= \EXP[ Y\,|\, \widehat{X}=\widehat{x}]$. Due to the identities 
\begin{align*}
    L^* - L^*(S) & = \Big( E[Y^2]- E\big[ m(X)^2\big ] \Big) - \Big(  E[Y^2]- E\big[ \widehat{m}(\widehat{X})^2\big ] \Big) \\
    &= \EXP\big[ (m(X)- \widehat{m}(\widehat{X}))^2 \big],
\end{align*}
the  regression null-hypothesis $\HH_0$  defined by (\ref{rtest1}) is equivalent to both 
\begin{align}
\label{rrtestl}
\EXP\big[ m(X)^2\big] =\EXP\big[ \widehat m(\widehat X)^2\big].
\end{align}
and 
\begin{align}
\label{aseq}
\PROB\big\{m(X)=\widehat m(\widehat X)\big\}=1.
\end{align}

As an estimate of $\EXP\big[ m(X)^2\big] -\EXP\big[ \widehat m(\widehat X)^2\big] $,
Devroye et al.  \cite[Section 3]{DeGyLuWa18} introduced the following  1-NN based test statistic
\begin{align}
\tilde T_n
&=\frac 1n \sum_{i=1}^n Y'_i (Y_{(n,1)}(X'_i)-\widehat Y_{(n,1)}(\widehat X'_i)).
\label{R3}
\end{align}
and accepted the  null-hypothesis $\HH_0$  if
\[
\tilde T_n\le \tilde a_n=\log n\left(n^{-1/2}+n^{-1/d}\right).
\]
For bounded  $Y$ and $X$, $X$ with  a density, and  $m$ satisfying the ordinary Lipschitz condition, \cite{DeGyLuWa18}  showed that this test is strongly consistent. For large $d$, the threshold above is too large to be of practical use.

\bigskip
We slightly modify the test statistic $\tilde T_n$ in a way that results in a negative bias. Define a nearest-neighbor based test statistic
\begin{align*}
T_n
&=
\frac 1n \sum_{i=1}^n \left(Y'_i  m_n(X'_i)-\widehat m_{n}(\widehat X'_i)^2\right),
\end{align*}
where  $ m_n$ and $\widehat{m}_n$  are the $k$-NN regression estimators defined in \eqref{est} and \eqref{est-1}, respectively. Again, here we set 
\begin{align*}
k &= 
k_n
=\lfloor \log n\rfloor,
\end{align*}
and with $a_n=\log n/\sqrt{n}$ as in (\ref{an}),
we  accept the regression null hypothesis if
\begin{align*}
T_n\le a_n.
\end{align*}

The following theorem is the main result of this section: 
\begin{theorem}
\label{expo*}
Let $d\ge 2$ and assume that ties occur with probability $0$. 

\begin{itemize}
\item[\rm (a)]
Under the regression alternative hypothesis
\begin{align*}
\lim_{n\to \infty} \PROB\{T_n\le a_n\}=0 .
\end{align*}
\item[\rm (b)]
If  $\widehat m$ satisfies the  modified Lipschitz condition,  $\EXP\big[ Y^4\,|\,  X\big] \le C$ a.s.\ for some  finite $C>0$, and the residual variance satisfies $\EXP\big[ (Y-\widehat m(\widehat X))^2\big] >0$, then
under the regression null hypothesis
\begin{align*}
\lim_{n\to \infty} \PROB\{T_n>a_n\}=0 .
\end{align*}
\end{itemize}
\end{theorem}

Similarly to Theorem~\ref{expo}, if in the proof of this theorem we assume that $Y$ is bounded and the Efron-Stein inequality is replaced by an exponential concentration inequality, then with the larger threshold $(\log n)^2/\sqrt{n}$ we obtain strong consistency.

\section{Proofs}

\label{sec-proofs} 
\subsection{Some moment inequalities}

A cone  with angle $\theta$ centered at the origin is the collection of all points $y\in \R^d$ such that $\text{angle}(y,z)\le \theta$ for some given $z
\in \R^d$. 
The following inequalities will be needed in the proofs of  Theorems~\ref{expo} and \ref{expo*}. 

\begin{lemma}
\label{ES}
Let $\gamma_d$ be the minimum number of cones centered at the origin and having angle  $\pi/6$ whose  union covers $\R^d$.
Set
\begin{align*}
m_n^*(x)
&=
\frac{1}{k_n}\sum_{i=2}^{n} Y_i\IND_{\left\{\mbox{$X_i$ is among the $k_n$ NNs of $x$ in $\{X_2, \dots ,X_n\}$}\right\}}~. 
\end{align*}
\begin{itemize}
\item[\rm (a)]
Under the assumptions of Theorem \ref{expo}, one has
\begin{align*}
\EXP\left[\left( \int |m_n(x)|\,\mu(dx)-\int |m^*_n(x)|\,\mu(dx)\right)^2\right]
&\le
\frac{16\gamma_d^2}{n^2}.
\end{align*}
\item[\rm (b)]
Under the assumptions of Theorem \ref{expo*}, one has
\begin{align*}
\EXP\left[\left(\int m(x)m_n(x)\, \mu(dx)-\int m(x)m^*_n(x)\, \mu(dx)\right)^2\right]
&\le
\frac{16C\gamma_d^2}{n^2}
\end{align*}
and
\begin{align*}
\EXP\left[ \left( \int m_n(x)^2\mu(dx)-\int m^*_n(x)^2\mu(dx)\right)^2\right] 
&\le
\frac{64C\gamma_d^2}{n^2}.
\end{align*}
\end{itemize}
\end{lemma}

\begin{proof}
We only prove the second half of Lemma \ref{ES}(b); the proofs of the other two inequalities are similar, but easier.
We have
\begin{align*}
\MoveEqLeft \EXP\left[ \left( \int m_n(x)^2\mu(dx)-\int m^*_n(x)^2\mu(dx)\right)^2\right] \\
&\le
\EXP\left[ \left( \int |m_n(x)+ m^*_n(x)|\cdot |m_n(x)- m^*_n(x)|\mu(dx)\right)^2\right] \\
&=
\EXP\Big[  \big|m_n(X_{n+1})+ m^*_n(X_{n+1})\big|\cdot \big|m_n(X_{n+1})- m^*_n(X_{n+1})\big|\\
&\quad \cdot \big|m_n(X_{n+2})+ m^*_n(X_{n+2})\big|\cdot |m_n(X_{n+2})- m^*_n(X_{n+2})\big|\Big],
\end{align*}
where $X_{n+1}$ and $X_{n+2}$ are  independent of the training samples $\D_n$ and have the common distribution of the $X_i$.  Therefore, using the notation in~\eqref{NN}, 
\begin{align*}
\MoveEqLeft
\EXP\left[ \left( \int m_n(x)^2\mu(dx)-\int m^*_n(x)^2\mu(dx)\right)^2\right] \\
&\le
\EXP\bigg[  \EXP\Big[  (|m_n(X_{n+1})|+ |m^*_n(X_{n+1})|)\frac{|Y_1|+|Y_{(n,k_n+1)}(X_{n+1})|}{k_n}\\
&\quad \cdot (|m_n(X_{n+2})|+ |m^*_n(X_{n+2})|)\frac{|Y_1|+|Y_{(n,k_n+1)}(X_{n+2})|}{k_n}\\
&\quad \cdot 
\IND_{\left\{\mbox{$X_1$ is among the $k_n$ NNs of $X_{n+1}$ and $X_{n+2}$ in 
$\{X_1, \dots ,X_n\}$}\right\}}\\
& \quad \quad \quad \Big| \,  X_1, \dots ,X_n,Y_1, \dots ,Y_n\Big] \bigg] \\
&\le
\frac{16Ck_n^2}{k_n^4} \\
& \quad \cdot \EXP\bigg[  \EXP\Big[ 
\IND_{\left\{\mbox{$X_1$ is among the $k_n$ NNs of $X_{n+1}$ and  $X_{n+2}$ in 
$\{X_1, \dots ,X_n\}$}\right\}}\\
&\qquad \quad \Big| X_1, \dots ,X_{n+2}\Big] \bigg] \\
&\le 
\frac{16C}{k_n^2}\\
&\quad \quad \cdot \PROB\left\{\mbox{$X_1$ is among the $k_n+2$ NNs of $X_{n+1}$ and  $X_{n+2}$ in 
$\{X_1, \dots ,X_{n+2}\}$}\right\}.
\end{align*}
Thus, 
\begin{align*}
\MoveEqLeft
\EXP\left[ \left( \int m_n(x)^2\mu(dx)-\int m^*_n(x)^2\mu(dx)\right)^2\right]\\
&\le
\frac{16C}{k_n^2(n+1)n}
\EXP\bigg[ \Big(\sum_{j=2}^{n+2}\\
&\qquad \IND_{\left\{\mbox{$X_1$ is among the $k_n+2$ NN's of $X_{j}$ in 
$\{X_1, \dots ,X_{n+2}\}$}\right\}}\Big)^2\bigg]\\
&=
\frac{16C}{k_n^2(n+1)n}
\EXP\bigg[ \Big(\sum_{j=2}^{n+2}\\
&\qquad \IND_{\left\{\mbox{$X_1$ is among the $k_n+1$ NN's of $X_{j}$ in 
$\{X_1, \dots ,X_{j-1},X_{j+1},\dots ,X_{n+2}\}$}\right\}}\Big)^2\bigg] \\
&\le
\frac{16C}{k_n^2n^2}((k_n+1)\gamma_d )^2,
\end{align*}
where for  the last step we refer to in \cite[Corollary 6.1]{GyKoKrWa02}.
\end{proof}

\begin{lemma}
\label{BLB}
(Extended Efron-Stein inequality \cite[Theorem 6]{BoLuBo04}.)
Let $A$ be a measurable set and $Z=(Z_1,\dots , Z_n)$ be an i.i.d.\ $n$-tuple of $A$-valued  random variables. 
Set $Z^{(1)}=(Z_2,\dots , Z_n)$.
Let $f$ and $g$ be real-valued measurable functions  on $A^n$ and $A^{n-1}$, respectively, such that   $f(Z)$  is integrable. Then,
\begin{align*}
\Var (f(Z))
&\le
n\EXP\big[ \, (f(Z)-g(Z^{(1)}))^2\, \big] .
\end{align*}
\end{lemma}

\subsection{Proof of Theorem~\ref{expo}}

\label{expo-proof} 
\begin{proof}[Proof of Theorem~\ref{expo}(a)]
Assume the classification alternative hypothesis and define $c^*$ as 
\begin{align}
\label{cstardef}
c^*=\EXP\left[ (|m(X)|-|\widehat m(\widehat X)|)\right]
&>
0.
\end{align}
For $n$ sufficiently large,
\begin{align*}
\PROB\{ T_n\le a_n\}
&\le
\PROB\{ T_n\le c^*/2\}
\le
\PROB\{ |T_n-c^*|\ge c^*/2\}
\le
\frac{2\EXP[\, |T_n-c^*|\, ]}{c^*}.
\end{align*}
Therefore, it suffices to show that 
\begin{align}
\label{i}
\lim_{n\to \infty} \EXP[\, |T_n-c^*|\,]  =  0.
\end{align}
We use the decomposition
\begin{align}
\label{18*}
T_n=(T_n -\EXP[\, T_n\mid \D_n\,] )+\EXP[ \, T_n\mid \D_n\, ].
\end{align}
One obtains
\begin{align}
\nonumber 
\EXP\big[ (T_n -\EXP[\, T_n\mid \D_n\,] )^2\mid \D_n \big]
&\le
\frac{2\EXP[Y^2]}{n}+\frac{2\EXP[m_n(X'_1)^2 \mid \D_n]}{n}\\
& \le 
\frac{4}{n}  \quad \text{a.s.}, \label{w2} 
\end{align}
and thus
\begin{align}
\label{ii}
\EXP\big[ \,  (T_n -\EXP[\, T_n\mid \D_n\,] )^2\big] 
&\to 0.
\end{align}
Furthermore, one has
\begin{align}
\label{TnDn}
\EXP[\, T_n\mid \D_n\,]
&=
\EXP\big[\,  Y'_1 \sgn(m_n(X'_1)) \mid \D_n\,\big] -\EXP\big[ \, |\widehat m_n(\widehat X'_1)| \,\big|\, \D_n\, \big]\nonumber\\
&=
\int m(x) \sgn(m_n( x))\, \mu(dx)-\int |\widehat m_n(\widehat x)|\,  \widehat \mu(d\widehat x).
\end{align}
We notice
\begin{align*}
| m(x) \sgn(m_n( x))-|m(x)|\, |
&\le
2| m_n( x)-m(x)|
\end{align*}
and obtain
\begin{align}
\label{iii}
\MoveEqLeft\EXP\left[ \; \Big|\int m(x) \sgn(m_n( x))\, \mu(dx)-\int |m(x)| \, \mu(dx) \Big|\; \right] \nonumber\\
&\le
2\EXP\left[ \int |m_n( x)-m(x)|\,  \mu(dx)\right]\nonumber\\
&\quad \to 0,
\end{align}
where the last step follows from  \cite[Theorem 6.1]{GyKoKrWa02}.
This theorem also yields
\begin{align}
\MoveEqLeft
\EXP\left[\,  \Big| \int |\widehat m_n(\widehat x)|\,  \widehat \mu(d\widehat x)-\int |\widehat m(\widehat x)| \, \widehat{\mu}(d\widehat x) \Big| \, \right] \nonumber \\
&\le
\EXP\left[ \; \int |\widehat m_n(\widehat x)-\widehat m(\widehat x)|\,  \widehat \mu(d\widehat x)\; \right] 
 \to 0. \label{iv}
\end{align}
From (\ref{TnDn}), (\ref{iii}),  (\ref{iv}),  and the definition of $c^*$ in \eqref{cstardef}  we get
\begin{align}
\label{v}
\EXP\big[ \,  \big|\EXP[ \, T_n\mid \D_n\, ] -c^*\big|\,\big] 
&\to 0.
\end{align}
Now (\ref{18*}), (\ref{ii}), and  (\ref{v}) yield (\ref{i})
and  part (a) of Theorem~\ref{expo} is proved.
\end{proof}

\begin{proof}[Proof of Lemma~\ref{MAD}]
With the notation $p=(a+1)/2$, one has 
\begin{align*}
\sum_{i=1}^n Z_i 
&=2B(n,p)-n,
\end{align*}
where $B(n,p)$ is a binomial random variable with parameters $n$ and $p\in [0,1]$.
Without loss of generality we assume that $\sigma^2=1-a^2=4p(1-p)>0$; otherwise (\ref{C4}) holds with equality. Therefore, $p\in (0,1)$.
Due to the chain of equalities 
\begin{align*}
\EXP \big|2B(n,p)-n\big|
&= \EXP \big|2(n-B(n,1-p))-n\big|\\*
&=\EXP \big|-2B(n,1-p))+n\big|\\
&=\EXP \big| 2B(n,1-p))-n\big| ,
\end{align*}
it suffices to consider the case of $p\le 1/2$, i.e., $a\le 0$.
We have 
\begin{align}
\EXP\left|\frac{1}{n}\sum_{i=1}^n Z_i \right|
&=
\EXP\left|\frac{2B(n,p)-2np}{n} +a\right|  \nonumber \\
&=
\EXP\left(\frac{2B(n,p)-2np}{n} +a\right)^+
+
\EXP\left(\frac{-2B(n,p) +2np}{n} -a\right)^+  \nonumber \\
&=
2\EXP\left(\frac{2B(n,p)-2np}{n} +a\right)^+
+
\EXP\left[ \frac{-2B(n,p)  +2np}{n} -a\right]   \nonumber  \\
&=
2\EXP\left[  \frac{ 2B(n,p)-n)^+}{n} \right]  +|a|    \label{lem-eq}   \\
&\ge 
2\PROB\big\{ B(n,p)/n> 1/2\big\} /n+|a|. \label{lem-low} 
\end{align}
Stirling's formula implies that
\begin{align*}
\PROB\big\{  B(n,p)/n> 1/2\big\}
&\ge
\frac{1}{\sqrt{2n}}e^{-nD\left( 1/2\Vert p\right) },
\end{align*}
where 
\[
D(\epsilon\Vert p)=\epsilon\ln\frac{\epsilon}{p}+(1-\epsilon)\ln\frac{1-\epsilon}{1-p}, 
\]
see  (4.7.2) on p.~115 in Ash \cite{Ash90}.
Since  $D\left( 1/2\Vert p\right)=-\ln \sigma$, we obtain 
\begin{align*}
\PROB\big\{  B(n,p)/n> 1/2\big\}
&\ge 
\frac{1}{\sqrt{2n}}\sigma^n,
\end{align*}
which, in view of \eqref{lem-low}, proves the lower bound in the lemma.

For the sake of simplicity, in the proof of the upper bound assume that $n$ is even.
As in the proof of the upper bound assuming that  $0<p< 1/2$,  we have 
\begin{align*}
\EXP\big[ ( B(n,p)-n/2)^+\big] 
&=
\sum_{n/2<j\le n}(j-n/2)\binom{n}{j}p^j(1-p)^{n-j}\\
&=
\sum_{0<j\le n/2}j\binom{n}{n/2+j}p^{n/2+j}(1-p)^{n/2-j}\\
&=
\sum_{0<j\le n/2}j\binom{n}{n/2+j}2^{-n}\left(\frac{p}{1-p}\right)^j 2^n\sqrt{p(1-p)}^n\\
&=
\sum_{0<j\le n/2}j\binom{n}{n/2+j}2^{-n}\left(\frac{p}{1-p}\right)^j \sigma^n.
\end{align*}
Therefore, in view of \eqref{lem-eq},   we have to prove that 
\begin{align*}
\sum_{0<j\le n/2}j\binom{n}{n/2+j}2^{-n}\left(\frac{p}{1-p}\right)^j 
&=
O\left(1/\sqrt{n}\right).
\end{align*}
Stirling's  formula implies that
\begin{align*}
\sqrt{\frac{n}{8j(n-j)}}
&\le 
\binom{n}{j}2^{-nh(j/n)}
\le
\sqrt{\frac{n}{2\pi j(n-j)}}\quad  \text{if \ $1\le j\le n-1$} 
\end{align*}
and
\begin{align*}
\binom{n}{j}2^{-nh(j/n)}
\le
1\quad  \text{if \ $0\le j\le n$},
\end{align*}
where, for $\epsilon\in (0,1)$, 
\begin{align*}
h(\epsilon)
&=
-\epsilon\log_2\epsilon -(1-\epsilon)\log_2(1-\epsilon) 
\end{align*}
is the binary entropy function,
see   p.~530 in Gallager \cite{Gal68}.
This implies 
\begin{align*}
\MoveEqLeft\sum_{0<j\le n/2}j\binom{n}{n/2+j}2^{-n}\left(\frac{p}{1-p}\right)^j\\ 
&\le 
\sum_{0<j\le n/4}j\sqrt{\frac{n}{2\pi (n^2/4-j^2)}}\left(\frac{p}{1-p}\right)^j \\
&\quad +
\sum_{n/4<j\le n/2}j2^{-n(1-h((n/2+j)/n))}\left(\frac{p}{1-p}\right)^j \\ 
&\le 
\sqrt{\frac{n}{2\pi (n^2/4-n^2/16)}}\sum_{0<j\le n/4}j\left(\frac{p}{1-p}\right)^j 
+
2^{-n(1-h(3/4))}n^2\\ 
&=
O\left(1/\sqrt{n}\right)+O\left(2^{-n(1-h(3/4))}n^2\right)\\
&=
O\left(1/\sqrt{n}\right),
\end{align*}
since  $p<1/2$.
\end{proof}

\begin{proof}[Proof of Theorem~\ref{expo}(b)]
Assume the classification null hypothesis.
We use the decomposition of $T_n$ in (\ref{18*}). The upper bound in 
(\ref{w2})  implies
\begin{align}
\label{w5}
 \PROB\big\{ T_n -\EXP[ \, T_n\mid \D_n\, ] \ge a_n/2 \big\}
&\le
\frac{8}{a_n^2 n}.
\end{align}
On the other hand, (\ref{TnDn}) and the classification null hypothesis (see \eqref{testl})  yield that
\begin{align}
\label{DD}
\EXP[\, T_n\mid \D_n\,] 
&\le
\int (|m(x)|-|\widehat m_n(\widehat x)| )\, \mu(dx)\nonumber\\
&=
\int (|\widehat m(\widehat x)|-|\widehat m_n(\widehat x)| )\, \widehat \mu(d\widehat x)\nonumber\\
&=
-\left(\int |\widehat m_n(\widehat x)| \, \widehat \mu(d\widehat x)
- \EXP\left[ \int |\widehat m_n(\widehat x)|\,  \widehat \mu(d\widehat x) \right] \right)\nonumber\\
&\quad +
\left( \int |\widehat m(\widehat x)| \, \widehat \mu(d\widehat x)
- \EXP\left[ \int |\widehat m_n(\widehat x)| \, \widehat \mu(d\widehat x) \right] \right).
\end{align}
Using first Lemma~\ref{BLB}  and the corresponding definition of~$\widehat m_n^*$,  and then part~(a) of Lemma~\ref{ES},  
one obtains
\begin{align}
\label{8'}
\Var\left( \int |\widehat m_n(\widehat x)|\, \widehat \mu(d\widehat x)\right)
&\le
n\EXP\left[ \left( \int |\widehat m_n(\widehat x)|\, \widehat \mu(d\widehat x)
-\int |\widehat m_n^*(\widehat x)|\widehat \mu(d\widehat x)\right)^2\right] \nonumber\\
&\le 
\frac{16\gamma_{d'}^2}{n}\nonumber\\
&=
\frac{C_d}{n}, 
\end{align}
where $C_d>0$ is finite.
We will  prove that
\begin{align}
\label{88'} 
\int |\widehat m(\widehat x)|\, \widehat \mu(d\widehat x)- \EXP\left[ \int |\widehat m_n(\widehat x)|\, \widehat \mu(d\widehat x)\right] 
<
0
\end{align}
for all  $n$  large enough.
We use the decomposition 
\begin{align}
\MoveEqLeft \int |\widehat m(\widehat x)|\, \widehat \mu(d\widehat x)- \EXP\left[ \int |\widehat m_n(\widehat x)|\, \widehat \mu(d\widehat x)\right] \nonumber \\
&=
\int |\widehat m(\widehat x)|\, \widehat \mu(d\widehat x)
-\int \big| \EXP[ \, \widehat m_n(\widehat x) \, ] \big|\, \widehat \mu(d\widehat x)
   \label{eq-dec1} \\
&\quad +\int \big|\EXP[ \, \widehat m_n(\widehat x) \,] \big|\, \widehat \mu(d\widehat x)
- \EXP\left[ \int |\widehat m_n(\widehat x)|\, \widehat \mu(d\widehat x)\right] .  \label{eq-dec2} 
\end{align}
For \eqref{eq-dec1}, under the modified Lipschitz condition the proof of Theorem 6 in  D\"oring et al.\ \cite{DoGyWa18} implies that
\begin{align}
\label{Do}
\int |\widehat m(\widehat x)|\, \widehat \mu(d\widehat x)
-\int \big|\EXP[\, \widehat m_n(\widehat x) \,] \big|\, \widehat \mu(d\widehat x)
&\le 
\int \big|\widehat m(\widehat x)-\EXP[ \, \widehat m_n(\widehat x)\,] \big|\, \widehat \mu(d\widehat x)\nonumber\\
&=O\left( \left(\frac{k_n}{n} \right)^{1/d'}\right).
\end{align}

Let  $R_{n,k_n}(\widehat{x})  =  \|\widehat x-\widehat X_{(n,k_n)}(\widehat x) \|$ be the $k_n$-NN distance and recall that ties occur with probability $0$ by assumption.   Since $\widehat{m}_n(\widehat x)$, defined in (\ref{est-1}), can be written in the form 
\[
\widehat{m}_n(\widehat x) = \frac{1}{k_n}\sum_{j=1}^{n} Y_j \, \IND_{\{\|\widehat x-\widehat X_j\|\le R_{n,k_n}(\widehat{x})  \}}\, , 
\]
one can see that given $R_{n,k_n}(\widehat{x})$,   the estimate $\widehat m_n(\widehat x)$ is conditionally distributed as the average of $k_n$  i.i.d.\ $\pm 1$-valued random variables. 
Conditioned on $ R_{n,k_n}(\widehat{x})$, these random variables have common variance 
\begin{equation}
k_n\cdot \Var\big( \widehat m_n(\widehat x)\mid  R_{n,k_n}(\widehat{x})  \big) 
=
1-\EXP\big[ \widehat m_n(\widehat x)\mid R_{n,k_n}(\widehat{x})  \big]^2.
\label{Vr}
\end{equation}
Then Jensen's inequality and Lemma \ref{MAD} imply 
\begin{align*}
\MoveEqLeft \int\EXP\big[ \, |\widehat m_n(\widehat x)|\, \big] \, \widehat \mu(d\widehat x)
 -\int\big| \EXP[\, \widehat m_n(\widehat x)\,] \big|\, \widehat \mu(d\widehat x)\\*
&\ge \int\EXP\left[ \EXP\big[ \, | \widehat m_n(\widehat x)| \, \big| \,   R_{n,k_n}(\widehat{x})  \big] \right] \, \widehat \mu(d\widehat x)\\
&\quad -\int\EXP\bigg[  \Big| \EXP[ \, \widehat m_n(\widehat x) \mid   R_{n,k_n}(\widehat{x}) \,  ] \Big| \bigg]  \, \widehat \mu(d\widehat x)\\
&\ge 
\frac{\sqrt{2}}{k_n^{3/2}}\int\EXP\left[ \left( 1-\EXP[\, \widehat m_n(\widehat x)\mid  R_{n,k_n}(\widehat{x}) \, ]^2\right)^{k_n/2}\right] \widehat \mu(d\widehat x).
\end{align*}
Again, apply Jensen's inequality twice:
\begin{align*}
\MoveEqLeft \int\EXP\big[ \ |\widehat m_n(\widehat x)|\, \big] \, \widehat \mu(d\widehat x)
 -\int\big| \EXP[ \, \widehat m_n(\widehat x)\, ] \big|\, \widehat \mu(d\widehat x)\\
&\ge 
\frac{\sqrt{2}}{k_n^{3/2}}\left(\int\EXP\Big[  1-\EXP[\, \widehat m_n(\widehat x)\mid  R_{n,k_n}(\widehat{x}) \, ] ^2\Big] \, \widehat \mu(d\widehat x)\right)^{k_n/2}\\
&\ge 
\frac{\sqrt{2}}{k_n^{3/2}}\left(\int\EXP\Big[  1-\EXP[ \, \widehat m_n(\widehat x)^2\mid  R_{n,k_n}(\widehat{x}) \, ] \Big] \, \widehat \mu(d\widehat x)\right)^{k_n/2}\\
&=
\frac{\sqrt{2}}{k_n^{3/2}}\left( 1-\int\EXP\big[ \, \widehat m_n(\widehat x)^2\big] \, \widehat \mu(d\widehat x)\right)^{k_n/2}\\ 
&\ge 
\frac{\sqrt{2}}{k_n^{3/2}}\left( 1-\int\widehat m(\widehat x)^2\widehat \mu(d\widehat x)+o(1)\right)^{k_n/2},
\end{align*}
where the last step holds because $\int\EXP\{\widehat m_n(\widehat x)^2\}\widehat \mu(d\widehat x)\to \int\widehat m(\widehat x)^2\widehat \mu(d\widehat x)$ as $n\to \infty$ by  \cite[Theorem 6.1]{GyKoKrWa02}.
By the condition $1-\int\widehat m(\widehat x)^2\widehat \mu(d\widehat x)= 1-\EXP[ \widehat m(\widehat X)^2] >0$, we 
therefore obtain 
\begin{align}
\label{S}
\int \big| \EXP[\, \widehat m_n(\widehat x) \,] \big|\, \widehat \mu(d\widehat x)
- \EXP\left[ \int |\widehat m_n(\widehat x)|\, \widehat \mu(d\widehat x)\right]
&\le 
-\frac{e^{-c_1k_n}}{k_n^{3/2}}
\end{align}
for some $c_1>0$. Since $k_n
=\lfloor \log n\rfloor$, 
(\ref{Do}) and (\ref{S}) yield
\begin{align*}
\MoveEqLeft \int |\widehat m(\widehat x)|\widehat \mu(d\widehat x)- \EXP\left[ \int |\widehat m_n(\widehat x)|\widehat \mu(d\widehat x)\right] \\
&\le 
O\left( \left(\frac{\sqrt{\log n}}{n} \right)^{1/d'}\right)-\frac{1}{n^{c_1/\sqrt{\log n}}(\log n)^{3/4}}\\
&<0
\end{align*}
if $n$ is large enough, and so (\ref{88'}) is verified.
From  (\ref{18*}),   (\ref{w5}), (\ref{DD}), (\ref{8'}) and (\ref{88'}), we now get for all $n$ large enough, 
\begin{align}
\PROB\left\{ T_n > a_n \right\}
&\le
\PROB\left\{ T_n -\EXP[\,  T_n\mid \D_n\, ] \ge a_n/2 \right\} \nonumber \\
&\quad +
\PROB\left\{  -\left( \int |\widehat m_n(\widehat x)|\widehat \mu(d\widehat x)-\EXP\left[ \int |\widehat m_n(\widehat x)|\widehat \mu(d\widehat x)\right] \right)> a_n/2 \right\} \nonumber \\
&\quad +\IND_{\{\int |\widehat m(\widehat x)|\widehat \mu(d\widehat x)- \EXP\left\{\int |\widehat m_n(\widehat x)|\widehat \mu(d\widehat x)\right\}>0\}}\nonumber \\
&\le
\frac{8}{a_n^2 n}
+\frac{4\Var\left( \int |\widehat m_n(\widehat x)|\widehat \mu(d\widehat x)\right)}{a_n^2} \nonumber \\
&\le
\frac{8 }{a_n^2 n}
+\frac{4C_d}{na_n^2} \label{eq-nod} \\
&\quad \to 0, \nonumber 
\end{align}
which yields part (b) of Theorem~\ref{expo}.
\end{proof}

\subsection{Proof of Theorem~\ref{expo*}}

\begin{proof}[Proof of Theorem~\ref{expo*}(a)]
Assume the regression alternative hypothesis.
We will prove that
\begin{align*}
T_n
&\to
\int m(x)^2 \mu(dx)-\int \widehat m(\widehat x)^2\, \widehat \mu(d\widehat x)
>0
\end{align*}
in probability.  We have 
\begin{align*}
T_n
&=
L_n-\widehat L_n,
\end{align*}
where
\begin{align*}
L_n
&=
\frac 1n \sum_{i=1}^n Y'_i  m_n(X'_i)
\end{align*}
and
\begin{align*}
\widehat L_n
&=
\frac 1n \sum_{i=1}^n \widehat m_{n}(\widehat X'_i)^2.
\end{align*}
Again, one can show that
\begin{align}
\label{ww}
\EXP\big[ \,  (T_n -\EXP[\, T_n\mid \D_n\,] )^2 \big] 
&\le
\frac{2\EXP[\,  {Y'_1}^2m_n(X'_1)^2\, ] }{n}+\frac{2\EXP[\,  \widehat m_n(\widehat X'_1)^4\, ] }{n}\nonumber\\
&=
\frac{4\EXP[\,  Y^4\, ](1+o(1)) }{n},
\end{align}
which yields 
\begin{align}
\label{tnconv-1} 
\lim_{n\to \infty} (T_n -\EXP[\, T_n\mid \D_n\,] )
&= 0
\end{align}
in probability.
One has
\begin{align*}
\EXP[ \, L_n\mid \D_n\,] 
&= \int  m( x)m_{n}( x) \, \mu(dx)
\to \int  m(x)^2\mu(dx)
\end{align*}
and
\begin{align*}
\EXP[ \, \widehat L_n\mid \D_n\, ] 
&= \int \widehat m_{n}(\widehat x)^2\, \widehat \mu(d\widehat x)
\to \int \widehat m(\widehat x)^2\, \widehat \mu(d\widehat x)
\end{align*}
in probability, because by  \cite[Theorem 6.1]{GyKoKrWa02}, 
\begin{align*}
 \int ( m_{n}( x)-  m( x))^2\, \mu(d x)
\to 0
\end{align*}
and
\begin{align*}
 \int (\widehat m_{n}(\widehat x)- \widehat m(\widehat x))^2\, \widehat \mu(d\widehat x)
\to 0
\end{align*}
in probability.  Thus we obtain that
\begin{equation}
\label{tnconv-2} 
\lim_{n\to  \infty} \EXP[ \, T_n\mid \D_n\,] =0
\end{equation}
in probability. 
Therefore, by \eqref{tnconv-1} and \eqref{tnconv-2},  under the alternative hypothesis
\begin{align*}
T_n
&\to
\int m(x)^2\mu(dx)-\int \widehat m(\widehat x)^2 \, \widehat \mu(d\widehat x)
>0\nonumber
\end{align*}
in probability, which proves  part~(a) of Theorem~\ref{expo*}.
\end{proof}

\begin{proof}[Proof of Theorem~\ref{expo*}(a)]
Assume the regression null hypothesis. We have 
\begin{align}
    \PROB\{T_n>  a_n\} & \le   \PROB\big\{ T_n -\EXP[\, T_n\mid \D_n\, ] \ge a_n/2\big\} \nonumber \\
    & \quad \mbox{}  + \PROB\big\{ \EXP[\, T_n\mid \D_n\, ] - \EXP[\, T_n\, ] \ge a_n/2\big\}   \nonumber  \\ 
    &\quad + \IND_{\{\EXP[\, T_n\, ]  > 0\} } \label{tnanbound} 
\end{align}

The bound 
(\ref{ww}) implies that
\begin{align}
\label{Tn}
 \PROB\big\{ T_n -\EXP[\, T_n\mid \D_n\, ] \ge a_n/2\big\}
&\le
\frac{8\EXP[\, Y^4\, ] }{a_n^2n}
\to 0.
\end{align}
Similarly to the proof of Theorem~\ref{expo}, using first Lemma~\ref{BLB}  and the corresponding definition of~$\widehat m_n^*$,  and then part~(b) of Lemma~\ref{ES},  
one obtains 
\begin{align}
\label{V1}
\MoveEqLeft \Var\left( \int  m( x)m_{n}( x) \, \mu(dx)\right)\nonumber\\
&\le
n\EXP\left[ \left( \int  m( x)m_{n}( x) \, \mu(dx)
-\int  m( x)m_{n}^*( x) \, \mu(dx)\right)^2\right] \nonumber\\
&\le 
\frac{16C\gamma_{d}^2}{n}
\end{align}
and
\begin{align}
\label{V2}
\Var\left( \int \widehat m_{n}(\widehat x)^2\, \widehat \mu(d\widehat x)\right)
&\le
n\EXP\left[ \left( \int \widehat m_{n}(\widehat x)^2\, \widehat \mu(d\widehat x)
-\int \widehat m_{n}^*(\widehat x)^2\, \widehat \mu(d\widehat x)\right)^2\right] \nonumber\\
&\le 
\frac{64C\gamma_{d'}^2}{n}.
\end{align}
Thus,
\begin{align}
\MoveEqLeft \Var\left(\EXP[\, T_n\mid \D_n\,] \right)\nonumber\\
&=
\Var\left( \int  m(x)m_{n}(x) \mu(dx)- \int \widehat m_{n}(\widehat x)^2\widehat \mu(d\widehat x)    \right)\nonumber\\
&\le
2\Var\left( \int  m(x)m_{n}(x) \mu(dx)\right)+ 2\Var\left(\int \widehat m_{n}(\widehat x)^2\widehat \mu(d\widehat x)    \right)\nonumber\\
&\le
\frac{c_d}{n} \nonumber 
\end{align}
with finite $c_d>0$. Therefore
\begin{align}
\label{wTn}
    \PROB\big\{ \EXP[\, T_n\mid \D_n\, ] - \EXP[\, T_n\, ]  
      \ge a_n/2\big\}   & \le \frac{4\Var\left(\EXP[\, T_n\mid \D_n\,] \right)}{a_n^2} \nonumber \\
      & \le \frac{4c_d}{n a_n^2}  \to 0. 
\end{align}

In view of  \eqref{tnanbound},  (\ref{Tn}),  and  (\ref{wTn}), it remains to prove that under the regression null hypothesis
\begin{align}
\label{Tp}
\EXP[\, T_n\, ] 
&< 0
\end{align}
if $n$ is large enough.
Under the  the null hypothesis (see \eqref{aseq}) one can use  the decomposition 
\begin{align}
\EXP[\, T_n\, ] 
&=
\int m(x) \EXP[\, m_{n}(x)\,] \, \mu(dx)-\int \EXP[\, \widehat m_{n}(\widehat x)^2] \,\widehat \mu(d\widehat x)\nonumber\\
&=
\int \widehat m(\widehat x) \EXP[\, m_{n}(x)\, ] \, \mu(dx)-\int \widehat m(\widehat x)^2\, \widehat\mu(d\widehat x)\nonumber\\
&\quad + \int \widehat m(\widehat x)^2\, \widehat \mu(d\widehat x)-\int \EXP[ \, \widehat m_{n}(\widehat x)\, ]^2 \, \widehat \mu(d\widehat x)\nonumber\\
& \quad -\int \Var\left(\widehat m_{n}(\widehat x)\right)\widehat \mu(d\widehat x).
\label{d32}
\end{align}
Again, under the modified Lipschitz condition the proof of Theorem 6 in  D\"oring et al.\  \cite{DoGyWa18} implies that
\begin{align}
\MoveEqLeft \int \widehat m(\widehat x) \EXP[ \, m_{n}(x)] \, \mu(dx)-\int \widehat m(\widehat x)^2 \, \widehat \mu(d\widehat x)\nonumber\\
&\le
\sqrt{\int \widehat m(\widehat x)^2\, \widehat \mu(d\widehat x)} \, 
\sqrt{\int\big(\EXP[\, m_{n}(x)] -\widehat m(\widehat x)\big)^2 \, \mu(dx)}\nonumber\\
&=
\sqrt{\int \widehat m(\widehat x)^2\, \widehat \mu(d\widehat x)} \, 
\sqrt{\int \big( \EXP[\, m_{n}(x)] - m( x)\big) ^2\mu(dx)}\nonumber\\
&=
O\left( \left(\frac{k_n}{n} \right)^{1/d}\right)
\label{d33}
\end{align}
and
\begin{align}
\MoveEqLeft \int \widehat m(\widehat x)^2\, \widehat \mu(d\widehat x)-
\int \EXP[\, \widehat m_{n}(\widehat x)]^2\, \widehat \mu(d\widehat x)\nonumber\\
&\le
\sqrt{\int \big(|\widehat m(\widehat x)|+|\EXP[\, \widehat m_{n}(\widehat x)]|  \big)^2 \, \widehat \mu(d\widehat x)}  \, 
\sqrt{\int \big(\EXP[\, \widehat m_{n}(\widehat x)] -\widehat m(\widehat x)\big)^2 \, \widehat \mu(d\widehat x)}\nonumber\\ 
&=
O\left( \left(\frac{k_n}{n} \right)^{1/d'}\right).
\label{d34}
\end{align}
Analogously to (\ref{Vr}) and by Jensen's inequality we obtain that
\begin{align*}
&k_n\cdot \int \Var\left(\widehat m_{n}(\widehat x)\right)\, \widehat \mu(d\widehat x)\\*
&\ge 
k_n\cdot \int\EXP\left[ \Var\left(\widehat m_{n}(\widehat x)\mid \widehat X_1,\dots ,\widehat X_n\right) \right] \, \widehat \mu(d\widehat x)\\
&=
\int\EXP\left[ \frac{1}{k_n}\sum_{j=1}^{k_n}\Var\left(Y_{n,j}(\widehat x)\mid \widehat X_1,\dots ,\widehat X_n\right) \right] \, \widehat \mu(d\widehat x)\\
&=
\int\frac{1}{k_n}\sum_{j=1}^{k_n}\EXP\left[ \left(Y_{n,j}(\widehat x)- \widehat m(\widehat X_{n,j}(\widehat x))\right)^2 \right] \, \widehat \mu(d\widehat x)\\
&\quad \to 
\EXP\big[ (Y-\widehat m(\widehat X))^2\big], 
\end{align*}
where for the limit relation we refer to Theorem 6.1 in \cite{GyKoKrWa02}.
Therefore, under the condition $\EXP\big[ (Y-\widehat m(\widehat X))^2\big] >0$, one obtains 
\begin{align}
 -\int \Var\left(\widehat m_{n}(\widehat x)\right)\widehat \mu(d\widehat x)
&\le 
-\frac{c_2}{k_n}
\label{d35}
\end{align}
with  $c_2>0$ for  $n$ large enough.
Thus, (\ref{d32}), (\ref{d33}), (\ref{d34}) and (\ref{d35}) yield 
\begin{align*}
\EXP\left[ T_n\right] 
&\le 
O\left( \left(\frac{1}{n} \right)^{1/d}\right)
+O\left( \left(\frac{k_n}{n} \right)^{1/d'}\right)
-\frac{c_2}{k_n}
< 0
\end{align*}
for $n$ sufficiently large (since $k_n
=\lfloor \log n\rfloor$) and so (\ref{Tp}) is verified. This completes the proof of part~(b) of Theorem~\ref{expo*}. 
\end{proof}

\section*{Declarations}

\noindent\textbf{Conflict of interest} The authors declare no conflict of interest.

\end{document}